\documentclass[runningheads]{llncs}
\usepackage[breaklinks, colorlinks, linkcolor=blue, citecolor=blue, urlcolor=black]{hyperref}
\usepackage{color}

\usepackage{graphicx}
\usepackage{tabularx}
\usepackage{tabularcalc}
\usepackage{booktabs}
\usepackage{multirow}
\usepackage[linesnumbered,ruled,vlined]{algorithm2e}

\SetCommentSty{AlgoCommentStyle}
\newcommand{\lIfElse}[3]{\lIf{#1}{#2 \textbf{else}~#3}}
\usepackage{algpseudocode}
\usepackage{xltabular}
\usepackage{amsmath}
\usepackage{amssymb}
\usepackage{xspace}
\usepackage{xcolor}
\usepackage{varioref}
\usepackage{cleveref}
\usepackage{aligned-overset}
\usepackage[shortlabels,inline]{enumitem}
\usepackage{boxedminipage}
\usepackage{comment}

\newcommand{\citep}{\cite}

\definecolor{pbblue}{RGB}{68,119,170}
\definecolor{pbcyan}{RGB}{102,204,238}
\definecolor{pbgreen}{RGB}{34,136,51}
\definecolor{pbred}{RGB}{238,102,119}

\usepackage{tikz}
\usetikzlibrary{arrows.meta}
\tikzset{%
  >={Latex[width=2mm,length=2mm]},
            base/.style = {rectangle, rounded corners, draw=black,
                           minimum width=4cm, minimum height=1cm,
                           text centered, font=\sffamily},
       startstop/.style = {base, fill=white,minimum width=4cm},
    activityRuns/.style = {base, fill=white,minimum width=4cm},
         process/.style = {base, minimum width=4cm, fill=white!15},
}

\newcommand{\solver}[1]{\textsc{#1}\xspace}
\newcommand{\technique}[1]{\textsc{#1}\xspace}
\newcommand{\scip}{\solver{SCIP}}

\newcommand{\vipr}{\solver{VIPR}}
\newcommand{\miplib}{\solver{MIPLIB}}
\newcommand{\soplex}{\solver{SoPlex}}
\newcommand{\bugger}{\solver{MIP-DD}}
\newcommand{\papilo}{\solver{PaPILO}}

\newcommand{\myorcidlink}[1]{\,\href{https://orcid.org/#1}{\raisebox{-0.45ex}{\includegraphics[width=1.8ex]{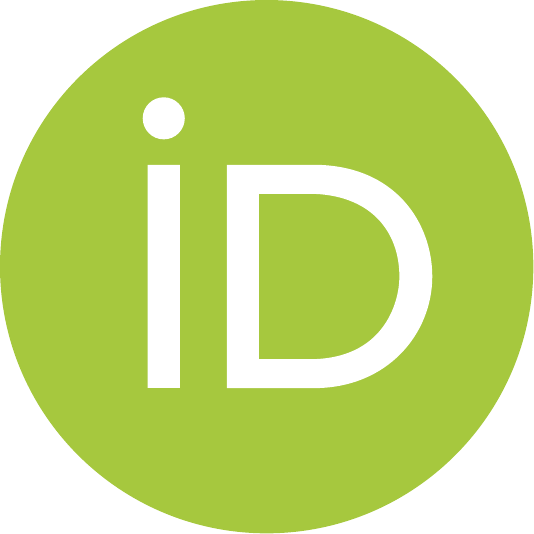}}}}

\usepackage{todonotes,xspace}

\makeatletter
\def\orcidID#1{\href{http://orcid.org/#1}{\protect\raisebox{-1.25pt}{\protect\includegraphics{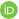}}}}
\makeatother

\begin{document}

\title{MIP-DD: Delta Debugging for Mixed Integer Programming Solvers}
\author{
    Alexander Hoen\inst{1}\orcidID{0000-0003-1065-1651}
    \and
	Dominik Kamp\inst{2}\orcidID{0009-0005-5577-9992}
    \and
    Ambros Gleixner\inst{1,3}\orcidID{0000-0003-0391-5903}
}
\institute{%
  Zuse Institute Berlin, Takustr.~7, 14195 Berlin, Germany\\
  \email{hoen@zib.de}
  \and
  University of Bayreuth, Universit\"atsstr.~30, 95440 Bayreuth, Germany\\
  \email{dominik.kamp@uni-bayreuth.de}\\
  \and 
  HTW Berlin, 10313 Berlin, Germany\\
  \email{gleixner@htw-berlin.de}
}

\authorrunning{A. Hoen, D. Kamp, A. Gleixner}
\maketitle

\begin{abstract}
  The recent performance improvements in mixed-integer programming (MIP) have been accompanied by a significantly increased complexity of the codes of MIP solvers, which poses challenges in fixing implementation errors.
In this paper, we introduce \bugger, 
a solver-independent tool, which to the best of our knowledge is the first open-source \emph{delta debugger} for MIP.
Delta debugging is a hypothesis-trial-result approach to isolate the cause of a solver failure.
\bugger simplifies MIP instances while maintaining the undesired behavior.
Preliminary versions already supported and motivated fixes for many bugs in the \scip releases 8.0.1 to 8.1.1.
In these versions, \bugger successfully contributed to 24 out of all 51 documented MIP-related bugfixes even for some long-known issues.
In selected case studies we highlight that instances triggering fundamental bugs in \scip can typically be reduced to a few variables and constraints in less than an hour.
This makes it significantly easier to manually trace and check the solution process on the resulting simplified instances.
A promising future application of \bugger is the analysis of performance bottlenecks, which could very well benefit from simple adversarial instances.

\end{abstract}
    
\section{Introduction}
\label{sec:intro}
Over the last few decades, solvers for \emph{mixed-integer programming}~(MIP) have witnessed remarkable performance improvements~\citep{KOCH2022100031}, allowing them to solve even large-scale problems despite the \textit{NP}-hardness of the problem class.
However, these performance gains come with increasingly complex algorithms.
For instance, in the last 10 years, the code base of \scip increased by approximately 44.6\%, precisely from 703\,965 lines of code in \scip~3.1 to 1\,017\,718 lines of code in \scip~9.1~\cite{SCIP9}.
This rise in complexity increasingly poses challenges in identifying and fixing implementation errors, so-called bugs.
Bugs may cause any kind of undesired solver behavior, especially claiming unboundedness or infeasibility falsely, or even returning ``optimal solutions'' that are actually infeasible or suboptimal~\citep{%
	CKSW13Hybrid,%
	Klotz2014,%
	PaxianBiere-POS23,%
	Steffy2011}.

In this paper, we describe a new open-source tool \bugger \citep{MIPDDCode} that can be used to facilitate and accelerate the process of locating implementation errors, ``bugs'', in any complete or heuristic MIP solver.
We discuss several successful case studies from the release process for Version~9.0 of the open-source MIP solver \scip~\citep{SCIP9}.

Let us first review common techniques used when debugging MIP solvers.\footnote{See also \url{https://scipopt.org/doc/html/DEBUG.php}.}
Assertions are boolean expressions that must evaluate to true and are used by programmers to test assumptions during runtime.
If an assertion fails, the run is immediately interrupted so that the invalid state can be inspected when the solver is compiled in debug mode and a software debugger is attached.
A debug solution mechanism allows the developer to enable further assertions which require a provided feasible solution to not be invalidated during the whole solving process.
Examples of such invalidations include pruning the node containing the debug solution or the application of a cut that renders the solution infeasible.

Despite their usefulness, these methods have their limitations.
First, due to computational overhead, assertions and debug solution checks are not present as frequently as necessary in order to pinpoint the exact moment an error occurs.
Second, bugs resulting from limited numerical precision might even pass assertions that incorporate certain tolerances.
Third, techniques in the MIP solvers, such as strong dual reductions~\citep{Achterberg07Thesis} or symmetry handling~\citep{Margot2010}, are allowed to cut off optimal solutions as long as at least one of them remains feasible.
For this reason, a failing debug solution check on some optimal solution does not automatically imply a solver error, which complicates efforts to investigate the issue of interest.

While experienced solver developers are trained to deal with these particular difficulties, tracing back the root cause of an error can still become costly or virtually impossible due to the sheer size and runtime of input problems encountered in many bug reports.
The goal of \bugger is to attack this size and time aspect of debugging by iteratively shrinking the problem size and the covered code base while preserving the reproducibility of the error.
The reductions performed by \bugger make it easier for a developer to isolate the actual bug on the reduced input and provide a targeted fix also for the original input.

In the literature, this approach is known as \emph{delta debugging}~\citep{ZNL99Deltadebugging},
which is an automated method to isolate the cause of a software failure, driven by a hy\-poth\-e\-sis-trial-result loop.
Delta debugging has been applied successfully
in the SAT and SMT solver community~\citep{BrunmayerBiere-FuzzingAndDeltadebuggingSMT-Solver,BrunmayerLonsingBiere10AutomatedTestingAndDebuggingSAT,KaufmannBiere-TAP22,NiemetzBiere-SMT13,PaxianBiere-POS23}.
To the best of our knowledge, this work describes the first application of delta debugging to mixed-integer programming and mathematical optimization in general.

The paper is organized as follows.
In \Cref{sec:description}, we describe the structure of our delta debugger.
In \Cref{sec:example}, we present a series of case studies to showcase the nature of different types of implementation errors and demonstrate the successful application of the delta debugger in these cases.
In \Cref{sec:conclusion}, we summarize different benefits of delta debugging and give an outlook on possible future improvements.

\section{Structure of the Delta Debugger}
\label{sec:description}

\begin{figure}[ht]
    \centering
    \scalebox{0.8}{
\begin{tikzpicture}[node distance=2.5cm,
    every node/.style={fill=white, font=\sffamily}, align=center]
    \node (start)[circle,fill=black,inner sep=3pt,minimum size=10pt,label=above:{Start}]{};
    \node (selectNext)[process, below of=start]{select modifier};
    \node (batches)[process, right of=selectNext, xshift=5.75cm, minimum width=9.5cm]{group admissible modifications into batches};
    \node (hasNext)[process, below of=batches, xshift=-2.75cm]{further batch?};
    \node (apply)[process, below of=hasNext]{apply batch};
    \node (callSolver)[process, below of=apply]{call solver};
    \node (evaluate)[process, right of=callSolver, xshift=3cm]{evaluate result};
    \node (reproduced)[process, above of=evaluate,] {bug detected?};
    \node (persist)[process, above of=reproduced]{revert batch};
    \node (terminate)[process, left of=hasNext, xshift=-3cm]{terminate?};    
    \node (stop)[circle,fill=black,inner sep=3pt,minimum size=10pt,label=below:{Stop},below of=terminate]{};
    \draw[->](start) -- (selectNext);
    \draw[->](selectNext) -- (batches);
    \draw[->](batches) -- ++(0,-1.25) -- ++(-2.75,0) -> (hasNext);
    \draw[->](hasNext) --node{Yes} (apply);
    \draw[->](hasNext) --node{No} (terminate);
    \draw[->](terminate) --node{No} (selectNext);
    \draw[->](apply) -- (callSolver);
    \draw[->](callSolver) -- (evaluate);
    \draw[->](evaluate) -- (reproduced);
    \draw[->](reproduced) --node{No} (persist);
    \draw[-](reproduced) -- ++(-2.75,0) --node{Yes}  ++(0,2.5) -- (hasNext.east);
    \draw[->](persist) -- (hasNext);
    \draw[->](terminate) --node{Yes} (stop);

    \draw[ultra thick] ([xshift=(-8),yshift=(8)]callSolver.north west) --node[above]{\textbf{Solver Interface}} ([xshift=(8),yshift=(8)]evaluate.north east) |- ([xshift=(8),yshift=(-8)]evaluate.south east) |- ([xshift=(-8),yshift=(-8)]callSolver.south west) -- cycle;
    
    \draw[ultra thick] ([xshift=(-8),yshift=(8)]batches.north west) --node[above]{\textbf{Modifier}} ([xshift=(8),yshift=(78)]persist.north east) |- ([xshift=(8),yshift=(-8)]reproduced.south east) |- ([xshift=(-8),yshift=(-8)]apply.south west) -- cycle;

    \draw[ultra thick] ([xshift=(-8),yshift=(8)]selectNext.north west) --node[above]{\textbf{Core Controller}} ([xshift=(8),yshift=(8)]selectNext.north east) |- ([xshift=(8),yshift=(-8)]terminate.south east) |- ([xshift=(-8),yshift=(-8)]terminate.south west) -- cycle;

  \end{tikzpicture}}
  \caption{Simplified workflow of \bugger. Termination criteria are shown in \Cref{algo:workflow}.}
\label{fig:simplified_workflow}
\end{figure}
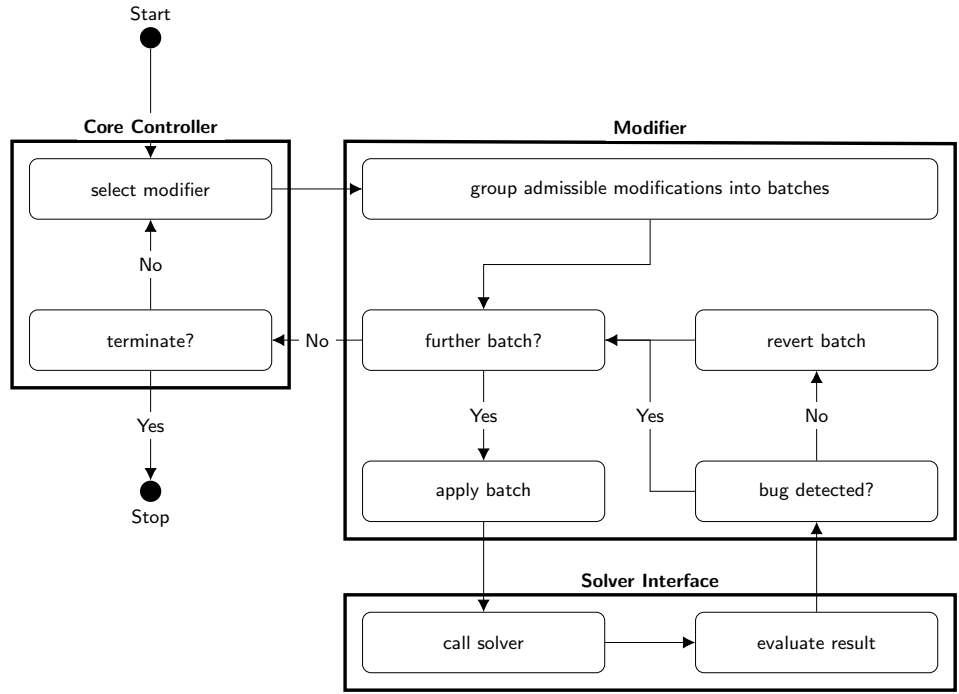

A simplified workflow of \bugger can be seen in \Cref{fig:simplified_workflow}.
\bugger consists of multiple \textit{Modifiers} each equipped with a unique strategy to isolate a bug.
Each modifier creates a local copy of the problem and modifies it.
In the general context of Delta Debugging, we then wish to test
the hypothesis that the solver exhibits the same bug also on the modified problem.
In order to accelerate this process, modifications can be grouped in batches and tested together as one hypothesis.
If the hypothesis is confirmed, i.e., the bug is reproduced, then all modifications are applied to the globally stored problem; otherwise, the modification is reverted. 

In the following subsections, we explain the implementation in further detail.
In \Cref{subsec::modifiers}, we first present the modifiers and their strategies.
In \Cref{subsec::underlying_solver}, we discuss the API by which the modifiers interact with the solver via a solver interface.
Through this API we invoke the solving of the (modified) problem and subsequent evaluation if a bug is reproduced such that a broad range of bugs can be detected, see also the example discussed in \Cref{subsec::count}.
In \Cref{subsec::workflow}, we motivate the ordering mechanism in which \bugger calls modifiers
to ensure that no potential modification was missed, in particular, if batches are used to speed up the process. 
Finally, in \Cref{subsec::recommendations}, we provide some general recommendations for using \bugger.

\subsection{Modifiers}
\label{subsec::modifiers}

\bugger consists of nine different modifiers each with a unique strategy to modify the problem or the settings in order to simplify the reproducing instance, reduce the code coverage, or accelerate the solving process to make the subsequent debugging easier.
Modifiers are selected sequentially in a fixed order.
To prioritize substantial reductions over more superficial reductions, the nine modifiers available in \bugger 1.0 are ordered as follows:

\begin{enumerate}
    \item Modifier \emph{Constraint} deletes a constraint from the problem.
    \item Modifier \emph{Variable} fixes a variable $x$ to the value of the reference solution. 
    \item Modifier \emph{Coefficient} deletes a coefficient $a_{ij}$ of a fixed variable and balances the left- and right-hand-side of the constraint $i$.
    \item Modifier \emph{Fixing} removes a fixed variable from the problem.
    \item Modifier \emph{Setting} switches a parameter of the solver to a user-defined value. These values are defined in the target setting file.
    \item Modifier \emph{Side} fixes inequalities $L \leq a \cdot x \leq U$ to $ax = a \cdot x^*$ with $x^*$ being the reference solution.
    \item Modifier \emph{Objective} sets the objective coefficient of a variable to zero.
    \item Modifier \emph{VarRound} rounds objective coefficient and the bounds of a variable to integer values to polish the reduced instance. Polishing facilitates debugging and rules out certain numerical difficulties.
    \item Modifier \emph{ConsRound} rounds the coefficients of a constraint to integer values to polish the reduced instance. The right-hand and left-hand sides are adjusted accordingly.
\end{enumerate}
Each modifier selects its modifications by iterating over the constraints and variables sequentially in the original order stored in the solver.

An important requirement is that each modification must preserve the feasibility of a provided reference solution.
The optimality of a reference solution is not required.
Therefore, in modifiers \emph{VarRound} and \emph{ConsRound}, the resulting bounds and sides are relaxed to include the reference solution if necessary.
Deleted variables are assumed to be fixed to the respective value of the reference solution.

Note that the algorithm behind \bugger is similar to the algorithm in~\cite[Algorithm 1]{Chinneck97FindingIIS} to compute an \emph{irreducible infeasible subset} (IIS)~\citep{Chinneck1991LocatingMI,Chinneck97FindingIIS,Pfetsch2002,GR90IIS}.
In the upcoming \bugger 2.1, an option to calculate an IIS based on the solver results is introduced.
When using this option, modifiers are restricted to relaxing problem reductions and only detecting infeasibility is considered as ``bug''.

By default, a solve call is triggered after each single modification.
Testing each modification separately, however, can be very time-consuming. 
To aggregate multiple modifications within a modifier the user can bound the number of solves per each invocation of the modifier by using the parameter \verb|nbatches|.
The modifiers then internally calculate how many modifications must be aggregated to a single hypothesis.
All modifiers apply a consecutive batch aggregation.
As an example consider the execution of the Constraint modifier for a MIP with $100\,000$~constraints.
With \verb|nbatches| set to 100, the modifier aggregates 1000 consecutive constraints to one batch and tests the removal of all these constraints at once, hence reducing the number of solve invocations to 100.

Because MIP solvers exhibit performance variability, most solvers provide the option to specify a random seed in order to trigger different solving paths.
Using this feature, one possible extension of \bugger could be to test each batch on multiple random seeds and accept the modification either if all of the seeds trigger the bug, or if at least one of the seeds triggers the bug.
The former variant could help to only accept modifications under which the bug is reproduced ``robustly''; the latter variant could in turn find more acceptable modifications, thereby yielding stronger reductions of the problem.

\subsection{Solver Interface}
\label{subsec::underlying_solver}

The communication with the solver happens in the solver interface, a general API that allows simple integration of any MIP solver.
The task of the solver interface is to report to the calling modifier whether a bug is reproduced by the current instance.
To do so, the solver interface applies several checks to detect the bugs of interest.
In the interface to \scip, the following checks are applied, which determine whether
\begin{itemize}
    \item the dual bound of the current problem cuts off the objective value of the reference solution. This is a contradiction to the feasibility of the reference solution. In this case, a suboptimality bug is detected.
    \item a solution returned by the solver is cut off by the problem restriction. In this case, an infeasibility bug is detected.
        For an unbounded problem, the solver may additionally provide a primal ray, which is then also verified to describe an unbounded improving direction.
    \item the primal bound is better than the evaluated objective value of the optimal solution, that is, the primal bound cuts off the best solution.
    \item an unexpected error occurs during the solving process, which is also considered a failure.
\end{itemize}
All these checks are provided by \bugger out of the box and can be used for any other MIP solver as well.

To integrate a solver into \bugger's API, a new class has to be created that inherits from a class \verb|SolverInterface|.
This solver interface interacts with the solver by
\begin{enumerate*}[label=(\alph*)]
    \item parsing the settings and instance files as well as loading the data into the internal data structure (\verb|parseSettings|, \verb|readInstance|),
    \item translating and loading the internal settings and problem to the solver (\verb|doSetup|),
    \item solving the problem and checking for bugs (\verb|solve|), and
    \item writing the internal settings and problem to files (\verb|writeInstance|).
\end{enumerate*}

The \verb|solve| method returns a \verb|signed char|, where the value 0 indicates that the solver finished and no bugs are detected.
Solver internal error codes must be mapped to negative values, whereas positive values are reserved for bugs detected by the solver-independent checks provided in the abstract SolverInterface class. 
\bugger provides template methods for checking the dual bound, the primal solutions, the primal ray, and the objective evaluation.
Unlike the debug solution assertions, the validity of these checks is not affected by strong dual reductions or symmetry handling.

\subsection{Workflow}
\label{subsec::workflow}

Aggregating modifications into batches can sometimes result in valid minor modifications being overlooked, especially when only a small subset of the elements within a batch is responsible for the solver's incorrect behavior.
Consequently, multiple calls to the same modifier may be required to uncover these hidden valid modifications.
This is particularly reasonable as the batch size typically decreases along with the problem size.
To address this, \bugger uses a stage-round mechanism to call the modifiers.
This approach is similar to the workflow of presolving, see, e.g.,~\cite{presolving_achterberg,GGHpapilo}.
In the $s$-th \emph{stage} only the modifiers with priority at most $s$ are called in sequence.
This is repeated over several rounds until no more changes are applied.
Then, the stage number is incremented or the algorithm terminates if there are no more modifiers to add to the stage cycle.
Hence, if in one round the problem remains the same, the next stage is entered to include an additional modifier.
Otherwise, the next round is entered to repeat the process with the same set of modifiers.

At the end of each round, the current settings-problem pair is written to files so that it can be used as an additional test instance to verify a bug fix.
The user can manipulate the mechanism by defining the initial and last stage as well as the maximal number of rounds.

\begin{algorithm}[ht]
	\DontPrintSemicolon
	\SetKwInOut{Input}{Input}\SetKwInOut{Output}{Output}
	\SetKwInOut{Init}{Initialization}
	\SetKwFunction{Goto}{go to}
	\SetKwFunction{And}{and}
	\Input{initial stage $s^0$, last stage $S$, last round $R$, a list of modifiers $M$ sorted by priority, a solver $\mathcal{M}$ which fails for the (settings, problem) pair $p^0$}
	\Output{set of reduced settings-problem pairs}
	$r \leftarrow 1$, $s \leftarrow s^0$\;
	\While{$r \leq R$ \And $s \leq S$}
	{
		$p^{r} \leftarrow p^{r-1}$\;
		\For{$t \leftarrow 1$ to $s$
                  }
		{
		  $p^{r} \leftarrow$ apply\_modifier($M_t, \mathcal{M}, p^{r}$)\;\label{algo:line:call_modul}
		}
                \lIfElse{$p^{r} = p^{r-1}$}
		{
            $s \leftarrow s + 1$
		}
        {
            $r \leftarrow r + 1$
        }
	}
    \Return{$\{p^{1}, \dots, p^{r}\}$}
	\caption{Core Controller of \bugger \label{algo:workflow}}
\end{algorithm}

\subsection{General Recommendations for Applying \bugger}
\label{subsec::recommendations}
\paragraph{Determine batch size.}
As explained in \Cref{sec:description}, the modifiers invoke the solver by default after each reduction.
Hence, a larger model size will increase the number of solve invocations.
For a more efficient process it is recommended to set parameter \verb|nbatches| which is maximum number of solves per modifier invocation.
This way, every successful solve invocation will simplify a certain chunk of the problem at once.

As a rule of thumb, in cases when reproducing a bug takes longer, parameter \verb|nbatches| should be set to a lower value.
This can help to reduce the total runtime of \bugger.
The final instance with a possibly reduced solving time can then be processed faster with a larger number of batches in order to find further reductions.

\paragraph{Define limits for solves.}

Some reductions can increase the solving time, especially when effective solving techniques are disabled.
Therefore, we recommend adding time or node limits for the underlying solver both to the original and target settings.
If no bug is detected within these limits, the corresponding reduction will be discarded immediately.
Typically, this speeds up the process and favors the creation of easy instances.

\paragraph{Suppress known issues.}

By default, all types of fails are treated as bugs.
However, it can happen that, for example, a dual fail turns into a primal fail during the reduction process.
Although discovering further issues is a beneficial side effect, it is sometimes desired to avoid particular types of fails in order to suppress unresolved known issues or prioritize certain bug fixes.
For this, the parameter \verb|passcodes| can be used.
All codes inside this list are additionally interpreted as correct results.
Negative values are reserved for solver-internal errors and positive values encode issues detected by \bugger.

\paragraph{Exact solution.}

For numerically sensitive issues it is especially important to ensure that the reference solution has feasibility violations as small as possible.
To illustrate potential confusion due to an inaccurate reference solution, we consider the simple example
\begin{align}
    \min_x\; -x_1&\nonumber\\
    x_1 + x_2 \leq &\,1 \\
    x \in \{0,1\}^2&\nonumber
\end{align}
with optimal solutions $x^* = (1,0)$.
We assume that the solver claims the suboptimal solution $x_f^* = (0,1)$ to be optimal.
Now, for the sake of illustration, assume that the infeasible solution $x^*_r = (1,1)$ would be used as a reference solution.  Then $x_2$ might be fixed to the reference value 1, which would render $x_f^*$ the correct optimal solution.
In this case, \bugger could interpret this correct result as a dual fail because the dual bound 0 claimed by the solver is higher than the objective value -1 of the reference solution assumed feasible.

Using the feasible solution $x^*=(1,0)$ as a reference solution instead, reliably avoids those problems.
Although slight feasibility violations of the reference solution are often compensated by the problem modifications, similar issues can still occur for an almost feasible reference solution due to differing scales when applying tolerances.
Therefore, using a feasible solution with minimal violation is recommended.
To generate an exactly feasible solution, an inexact solution can be either polished or an exact optimal solution can be computed directly, for example by using the numerically exact version of \scip~\citep{EiflerGleixner2022_Acomputational,EiflerGleixner2023_Safeandverified,exscipgithub}.

\section{Case Study}
\label{sec:example}
The primary objective of Delta debugging is to enable rapid bug tracking and save time during the debugging process.
While these time savings are difficult to quantify directly, we use five case studies in this section to demonstrate that \bugger often produces significant reductions leading to an effective acceleration and motivation in the debugging process.
During the development of the \scip versions 8.0.1 to 8.1.1, preliminary versions of \bugger were applied to significantly reduce the reproducing instances of 24 among all the 51 MIP-related and finally fixed bugs listed in the release reports.
Further statistics on the bugs can be found in the repository~\citep{MIPDDCode}.
Moreover, reduced instances generated by these bugs have been added to the set of instances that are tested regularly for sustainability during \scip development.
They can be found in the \scip repository\footnote{\url{https://github.com/scipopt/scip/tree/master/check/instances/Issue}}.

Before presenting specific examples, we need to briefly discuss the fact that MIP solvers typically employ floating-point arithmetic.
In floating-point arithmetic, a limited amount of memory is allocated to represent a number.
Hence, not all numbers can be represented accurately.
If more digits are needed than allocated the last digits are omitted leading to a loss in precision.
Among further problems, calculation in floating-point arithmetic is neither associative nor distributive meaning that $(a+ b) + c$ is not necessarily equal to $a + ( b + c )$ and $ (a + b ) \cdot c$ may yield a different result than $ a \cdot c + b \cdot c$.
For further information about floating-point arithmetic, we refer to \cite{Goldberg91_WhatEveryComputerScientistShouldKnow}.

To deal with these issues, floating-point MIP solvers typically define a \emph{zero tolerance}~$\varepsilon$  and a \emph{feasibility tolerance}~$\delta$ to overcome such shortcomings.
If two numbers $a,b$ are within $\varepsilon$-range, i.e., $|a-b|<\varepsilon$, they are considered equal.
In \scip, a constraint $a^T x \leq b$ with a solution $x^*$ is considered feasible if
\[
\dfrac{a^T x^* - b}{\max\{1, |b|, |a^T x^*|\}} < \delta.
\]
This complicates the development of a MIP solver and increases the potential for incorrect results on numerically difficult instances.



In the following case studies, we use \bugger 1.0 (Hash 455b7913) 
on a standard MacBook Air with an Apple(R) Silicon(TM) M2 2022 8‑Core CPU @ 3.49 GHz, 8‑Core GPU, and 16 GB RAM.
All times are given in seconds.


\begin{table}\small
    \begin{tabular*}{\textwidth}{@{}l@{\;\;\extracolsep{\fill}}lrrrrrrrrr}
        \toprule
        Case
        & \multicolumn{3}{c}{Original Instance}
        & \multicolumn{3}{c}{Final Instance}
        & \multicolumn{3}{c}{MIP-DD Statistics}\\
        \cmidrule{2-4} \cmidrule{5-7} \cmidrule{8-10}
        & Vars & Conss & Nonzeroes 
        & Vars & Conss & Nonzeroes
        & Rounds & MIP Solves & Time [s]\\
        \midrule
        \ref{subsec::dualinfer}
        & 6 & 6 & 36
        & 5 & 3 & 15
        & 12 & 152 & 0.2\\
        \ref{subsec::3499_rococo}
        & 3117 & 1293 & 11751
        & 7 & 4 & 11
        & 66 & 69430 & 2869 \\
        \ref{subsec::3662}
        & 653 & 745 & 1819 
        & 7 & 2 & 8
        & 10 & 2318 &4.1\\
        \ref{subsec::normalize} 
        & 653 & 745 & 1819 
        & 4 & 4 & 11 
        & 8 & 2198 & 4.3\\
        \ref{subsec::count}
        & 204 & 104 & 940
        & 0 & 0 & 0
        & 5 & 638 & 3.4\\
        \bottomrule
    \end{tabular*}
    \caption{Problem size reductions and MIP-DD performance.}
    \label{tab::stats}
\end{table}

\subsection{Modifying the Settings to Detect Wrong Calculation of Activities in DualInfer}
\label{subsec::dualinfer}

We first want to highlight the ability of \bugger to identify faulty components by manipulating the settings.
For this, we take a look at an instance reported by a \scip user that has been in the backlog for over three years (\scip git hash \verb|d4c0ad7644|).
Using the so-called aggressive presolve settings leads to an incorrect solution while default settings result in the correct solution. 
To identify the faulty presolve routine we define the default settings as the target setting and run \bugger.
As a result, all settings are reset to their default except for \technique{Dual-Infer}~\citep{presolving_achterberg}, which remained activated pointing to this presolver as the faulty component.

In fact, the reason for the incorrect solution is
that \technique{Dual-Infer} incorrectly calculated the min-max-residuum for the subset $\hat{N_k} = N\backslash\{k\}$ where $N$ represents the index set of all variables and $k \in N$.
If applied correctly, for a constraint $\sum_{i\in N} a_{i} x_i \circ b$ with $\circ \in\{=,\leq,\geq\}$ the maximal activity $a_{max}^k$ is calculated by
\begin{align*}
    a_{max}^k = \sum_{i\in \hat{N_k}, a_{i}>0} a_{i}\cdot u_i + \sum_{i\in \hat{N_k}, a_{i}<0} a_{i}\cdot l_i \;.
\end{align*}
In the reduced instance
\begin{align*}
\begin{aligned}
  \min_{x\geq0} 25\, x_0&\\
\begin{array}{rrrrr}
    51.2 x_{0}       & + 25.6 x_{1} & +44.8\, x_2& +38.4 \, x_3  & + 44.8 \, x_4 \\
    4 \, x_0 & +5 \, x_1 & -8 \, x_2  & -x_3 &-2\, x_4 \\
     -1.024\, x_0 &+2.688\, x_1 &-0.896\, x_2 &-0.768\, x_3 &+5.504\, x_4  \\
\end{array}&
\begin{array}{l l}
    = & 192000 \\
    \leq & 0 \\
    \leq & 0 \\
\end{array}\\
     x_0\leq 1000,\; x_1 \leq 3000,\; x_2 \leq 2000,\; x_3 \leq 800,\; x_4 &\leq 2600 
  \end{aligned}
\end{align*}
the maximum residuum $a_{max}^0$ for the first constraint is $313\,600$.
However, \scip claims that the residuum is $153\,600$ leading to incorrect assumptions in the subsequent solving process.
Since \technique{Dual-Infer} applies dual reductions using the debug solution is unlikely to be effective in locating the bug even if called frequently during presolving. 
The reason is an index shift where instead of $u_i$ and $l_i$, $u_0$ and $l_0$ were always applied. 
This is fixed with commit \verb|cf362110|.

\subsection{Fixing Invalid Separation of Interior Solutions}
\label{subsec::3499_rococo}

In this section, we want to analyze a bug that appeared on the \miplib\citep{MIPLIB2017} instance \verb|rococoC10-001000| when enabling the non-default separators intobj, closecuts, and oddcycle in the root node.
This bug can be replicated with \scip commit \verb|82835b0d|.

For this problem with a known true optimal solution value of 11460, a solution with value 23706 was declared optimal after 5776 LP iterations.
Applying delta debugging with disabled separators, heuristics, and presolvers as target settings, as well as suppressed primal fails, and 1,000 batches, led to the small integer program

\begin{align*}
    \begin{aligned}
        \min_{x \in \mathbb{N}_0} x_{7377}\\
        -27856 x_{6324} - 34000 x_{6326} +x_{6331} &= 6144\\
        -1289 x_{6324} -1916 x_{6326} + x_{6333} &= 627\\
        - x_{6333} + x_{7377} &= 9544 \\
        -1346 x_{3740} - 426 x_{5494} + x_{6331} &\geq 32642 \\
        x_{3740},x_{5494},x_{6324}, x_{6326}&\leq 1\\
        x_{6333}\leq 5748,\;x_{6331} \leq 102000,\;x_{7377} &\leq 262689\\
    \end{aligned}
\end{align*}

\noindent
in less than an hour.
The small program has the same optimal value as the original problem yet \scip returns an ``optimal'' solution of value 12087 after 4 LP iterations.
The only remaining non-disabling settings are presolving/maxrounds = -1, constraints/linear/max\-pre\-rounds = -1, separating/close\-cuts/freq = 0, and separating/intobj/freq = 0.
This means that apart from some presolving only separators closecuts and intobj are required to reproduce this issue.
As closecuts is known as a meta separator, these settings indicated that there might be an issue in separator intobj when called by separator closecuts~\citep{Gomory63AlgorithmIntegerSolutions,Gomory58CuttingPlane,cip_phd}.
Separator closecuts determined an interior point of the LP relaxation which is handed to intobj to be separated.
Separator intobj relies on the integrality of the objective function, which is exploited by cutting off a solution with a fractional objective value.
This is done by rounding this value up and using it as left-hand side of the objective row.
The approach implicitly assumes that there is no feasible solution with a smaller objective value than the infeasible solution to be separated since all undiscovered better solutions would otherwise be cut off.
However, for the arbitrary LP-interior points computed by separator closecuts, this does not hold in general.
It turned out that the value of the separated solution is approximately 11646.03 indeed exceeding the optimal value 11460.
Restricting intobj to optimal solutions of LP relaxations finally resolved this issue.

During preliminary experiments, we observed on this issue that using default settings as target settings, results in a larger final instance with 1528 variables and 301 constraints.
The runtime was even comparable to the previously presented setup although disabling plugins usually leads to larger solving times of the single solves.
An explanation for this behaviour might be given by the complexity of the solver under default settings which under certain circumstances requires a more complex external problem to reproduce the failing internal state.
Nevertheless, targeting default settings is desirable for the generation of lightweight but flexible regular tests.
And for other issues, comparably small reproducing instances could be found.


\subsection{Detecting Invalid Propagations in the Presolving Library \papilo}
\label{subsec::3662}
It is worth mentioning that \bugger can also reveal further bugs that are not even reproducible initially, especially in instances that are numerically challenging in general.
\papilo~\citep{GGHpapilo}, a presolving library for integer and linear optimization, which is part of \scip's presolving routines, was not called at all for the instance reported by a user.
However, the delta-debugger reduced it to 
\begin{align*}
    \begin{aligned}
        -1.64216813092803 \,\cdot 10^{-8}\, x_{350} + x_{351} &= 0\\
        0.5 x_{179} + 4 x_{183} + 0.5 x_{189} +1.5 x_{195} +6.3 x_{240} +x_{351} &= 6.29999999999997\\
        x_{351} \geq 0\;, x_{350}, x_{183} , x_{189}, x_{240}, x_{179} &\in \{0,1\},
    \end{aligned}
\end{align*}
where \papilo (githash \verb|df30ae49|) states infeasibility although a solution within $\delta$-tolerance exists, namely $x_{240} = 1$ and otherwise zero.

\papilo chooses to apply \emph{Probing}~\citep{Savelsbergh_Preprocessing_and_Probing_Techniques,presolving_achterberg} to the variable $x_{240}$.
\technique{Probing} tentatively fixes binary variables to 0 and 1, applies constraint propagation~\citep{Savelsbergh_Preprocessing_and_Probing_Techniques} and tries to conclude implication based on the outcomes of the propagation.
By fixing $x_{240}$ to 0, the problem correctly propagates to infeasibility.
If the variable $x_{240}$ is fixed to $1$, propagation tries to find new bounds for the variable $x_{350}$ by resolving the constraints by $x_{350}$:
\begin{align*}
    x_{350} \leq\Bigg\lfloor \frac{4.59999999999998-4.6}{1.64216813092803 \,\cdot 10^{-8}}\Bigg\rfloor\; .
\end{align*}
Since the quotient exceeds the $\varepsilon$ tolerance of $1e-9$ it is floored leading to the new upper bound of -1.
Since -1 is smaller than the lower bound, fixing $x_{350}$ leads to infeasibility also for $x_{240} = 1$.
Combined with the infeasibility of $x_{240} = 0$ \technique{Probing} declares global infeasibility.

To resolve this undesired numerical effect, flooring and ceiling for integer variables are applied with more caution.
In this case, $\frac{4.59999999999998-4.6}{1.64216813092803 \, 10^{-8}}$ is first ceiled to 0 against the feasibility tolerance $\delta$.
Then it is checked if $x_{350}=0$ still satisfies the original constraint within $\delta$.
If this is the case, this ceiled value is considered the upper bound.
Otherwise, we can safely cut off this value and the floored value is a valid upper bound.
This was fixed with the commit \verb|019915e0|.

\subsection{Identifying a Numerically Critical Constraint Normalization}
\label{subsec::normalize}

Since \papilo is not called in the original instance fixing the bug in \papilo does not fix the incorrect solver status of the original instance discussed in \Cref{subsec::3662}.
Therefore, re-running \bugger with the fixed \papilo and suppressed primal fails results in a new reduced problem 
\begin{align*}
    \begin{aligned}
        \begin{array}{rrrr}
            -6.88914620404344 \, x_{383} &+ x_{384} &&\\
            -20000007 \, x_{383} & + x_{384} & + x_{385} & + x_{386} \\
             & x_{384} &- x_{385} &+ x_{386}\\
             &&& x_{386}\\
        \end{array}&
        \begin{array}{ll}
             = & 0 \\
             \leq & 0\\
             = & 0 \\
             = & 0.599028894874692
        \end{array}\\
         x_{383} \in \{0,1\}, x \geq 0 &,
    \end{aligned}
\end{align*}

\noindent
where \scip (git hash \verb|cf244e14|) returns a false solution.

In this problem, the second constraint leads to numerical difficulties.
During presolving, all variables except $x_{383}$ are either fixed ($x_{386}= 0.599028894874692$ and $x_{385} = x_{386}$) or substituted resulting in the following problem 
\begin{align}
    \label{eq::3662::normalized_2}
     -20000007\, x_{383} \leq -1.198057789749384 \\
     x_{383} \in \{0,1\}\, .\nonumber
\end{align}
\scip normalizes the remaining constraint by dividing it by $-20000007$ leading to a right hand side  that is smaller than the feasibility tolerance $\delta$:
\begin{align}
    \label{eq::3662::normalized_3}
     x_{383} \geq 5.990286852146522 \cdot 10^{-8} \\
     x_{383} \in \{0,1\}\, .\nonumber
\end{align}
Hence, \eqref{eq::3662::normalized_3} implies the lower bound of $x_{383}$ and can be deleted due to the redundancy allowing to choose any value for $x_{383}$ within its domain.
Setting $x_{383}$ to 0 invalidates \eqref{eq::3662::normalized_2} since $0 + 1.198057789749384\not \leq 0$. 
If the normalization for vanishing left-hand sides is skipped, coefficient tightening~\citep{Savelsbergh_Preprocessing_and_Probing_Techniques,presolving_achterberg} for \eqref{eq::3662::normalized_2} is applied leading to 
\begin{align}
    \label{eq::3662::normalized_4}
     -1.198057789749384 \, x_{383} \leq -1.198057789749384\\
     x_{383} \in \{0,1\}\, .\nonumber
\end{align}
Immediately after the tightening, normalization is invoked again, normalizing \eqref{eq::3662::normalized_4} to $x_{383} \geq 1$ which implies the correct fixing of $x_{383}$ to 1.
This is fixed with the commit \verb|6b2dd64c|.


\subsection{Fixing the Solution Count}
\label{subsec::count}

\bugger can not only be used to detect primal or dual fails but can serve for testing the entire solver functionality.
For example, \scip incorporates a counting mode, in which the number of all feasible value assignments on integral variables is counted\footnote{\url{https://www.scipopt.org/doc/html/COUNTER.php}}.
If a problem is already solved during presolving, a solution count of zero was claimed, even though a feasible solution was provided when using the regular optimization mode.
This bug can be reproduced with \scip commit \verb|ca23f99f|.
We extended the \scip interface to return an error if \scip claims a solution count of zero on a feasible instance.
\bugger showed that this issue already occurs for the trivial problem with neither variables nor constraints due to the mishandling of a presolved problem in which all variables are fixed.
With this information the issue could be resolved with the commit \verb|434945be|.

\section{Conclusion and Outlook}
\label{sec:conclusion}

To conclude, let us revisit the different benefits of delta debugging during the development process of an MIP solver and give an outlook on possible enhancements for the future.
%
First, as shown in \Cref{sec:example}, \bugger is able to fulfill its general purpose of reducing the instance size while still providing incorrect output.
This spares a lot of time during debugging by reducing the hurdles for investigating issues because all relevant calculations done by the solver on the reduced instance can usually be traced and checked manually.
Further, \bugger can be applied even if other techniques such as providing a debug solution fail.

Second, in the context of collaborations with industry, it regularly happens that instances cannot be shared due to the General Data Protection Regulation~\citep{GDPR2016a} or other confidentiality reasons.  In this case, a reduced instance with renamed variable and constraint names obtained by delta debugging may barely contain any sensitive information anymore about the original instances, making the reduced instance safe for customers to share with developers.

Third, \bugger itself does not prescribe what counts as an error. Instead, the definition of an error is provided by the implementation of the solver API.
Hence, \bugger can not only be applied to detect bugs caused by primal or dual infeasibility but can analyze any unwanted behavior specific to a certain solver.
As one example, the \bugger 2.0 will contain, along with improvements for automatic limit settings, adaptive batch numbers, and extended arithmetical options, an experimental interface to the numerically exact version of  \scip~\citep{EiflerGleixner2023_Acomputational,EiflerGleixner2024_Safeandverified,exscipgithub} and the verification software \vipr for checking validity of certificates~\citep{CheungGleixnerSteffy2017}.
This will make it possible to identify the generation of an invalid certificate as another type of failure.
In this context, reducing the size of problem instances is helpful because it typically goes hand in hand with reducing the size of invalid certificates, which therefore may become easier to analyze.

In an even broader context, \bugger could also be utilized to track down performance bottlenecks, e.g. by letting the solver return an error at a certain tree depth which should not be reached on a given instance. This way, issues with exploding tree sizes on instances turning out to be simple could be investigated manually, which contributes to creating more performant solvers in itself.

\bigskip

\noindent
\textbf{Acknowledgements.}
The work for this article has been partly conducted within the Research Campus Modal funded by the German Federal Ministry of Education and Research (BMBF grant numbers 05M20ZBM) and the Bayreuth Research Center for Modeling and Simulation (MODUS) funded by the German Research Foundation (DFG project RA 1033/3-1).
\\
We would like to thank all who contribute to \scip by submitting bug tickets and allowing us to publish their instances.
Furthermore, we wish to thank Leon Eifler for helping to implement the interface to the exact version of \scip and \soplex, Mark Turner and Felipe Serrano for their careful proofreading, and Armin Biere for sharing his experience and perspective on delta debugging in the SAT solver world. 
Finally, we would like to thank the associate editor and two anonymous reviewers for their insightful comments and suggestions, which have contributed to improving the quality of this work.

\newpage

\bibliographystyle{unsrt}
\bibliography{bibliography,refArticles,refBooks,refOther}

\end{document}